\documentclass{article} 
\usepackage{graphicx}
\usepackage{latexsym}
\usepackage{amssymb}
  \newtheorem {korolar}    {Corollary}  
 \newtheorem {propozicija}{Proposition}

  \newcommand {\dzcite} {\mbox{\bf Proof : }}

\newcommand{\qed}{$\square$}
  
  \newcommand {\tj} {\'{c}}

  \newcommand {\ra} {\rightarrow}
  \newcommand {\lra} {\leftrightarrow}
  \newcommand {\longra} {\longrightarrow}

\title{Logic of Paradoxes in Classical Set Theories}

\author{Boris Čulina  \\  Department of Mathematics \\ University of Applied Sciences Velika Gorica
	\\ Zagreba\v{c}ka cesta 5, 10410 Velika Gorica, Croatia
	\\ e-mail: boris.culina@vvg.hr}
\date{}

\begin {document}
  \maketitle
\textit{Logistic is not sterile; it engenders antinomies.}

\rightline{H. Poincar\'e}
\begin{abstract}

    According to Georg Cantor \cite{can1,can} a set is any multitude
 which can be thought of as  one ("jedes Viele, welches sich als Eines denken l\"{a}{\ss}t")
 without contradiction -- a consistent multitude. Other multitudes
 are inconsistent or paradoxical. Set theoretical paradoxes have
 common root --  lack of understanding why some multitudes are
 not sets. Why some multitudes of objects of thought cannot
 themselves be objects of thought? Moreover, it is a logical truth that such multitudes do exist. However  we do not understand this logical truth so well as we understand, for example, the logical truth $\forall x \ x=x$. In this paper we formulate a logical truth which we call the productivity principle. Bertrand Rusell \cite{rus} was the first one to formulate this principle, but in a restricted form and with a different purpose. The principle  explicates a logical mechanism that lies behind paradoxical multitudes, and  is understandable as well as any simple logical truth. However,  it does not explain the
  concept of set. It only sets  logical bounds of the concept within the framework of the classical two valued $\in$ - language. The principle behaves as a logical regulator of  any  theory we formulate to explain and describe sets. It provides tools to identify paradoxical classes inside the theory. We  show how the known paradoxical classes
follow from the productivity principle and how the principle gives us a uniform way
to generate new paradoxical classes. In the case of $ZFC$ set theory the productivity principle shows  that the limitation of size principles are of a restrictive nature and that they do not explain which classes are sets.   The productivity principle, as a
  logical regulator, can have a definite heuristic
  role in the development of a consistent set theory. We sketch such a theory -- the cumulative cardinal theory of sets. The theory  is based
  on the idea of cardinality of collecting
  objects into sets. Its development is guided by means of the productivity principle in such a way that  its consistency seems plausible.  Moreover, the theory inherits good properties from cardinal conception and from cumulative conception of sets. Because of the cardinality principle it can easily justify the replacement axiom, and because of the cumulative property it can easily justify  the power set axiom and the union axiom. It would be possible to prove that the cumulative cardinal theory of sets is equivalent to the Morse -- Kelley set theory. In this way we provide a natural and plausibly consistent axiomatization for the Morse -- Kelley set theory.

\end{abstract}
  \bibliographystyle {alpha}

 \section{The productivity principle}

In this section we will formulate a simple logical criterion to distinguish between sets and proper or paradoxical classes. We have named it the productivity principle.

When we talk about some objects, it is natural to talk about
collections of these objects as well. The same happens when we talk
about collections themselves. But the moral of paradoxes is that we
can not talk freely about their collections.

First of all, the language has to be made precise. Since we are talking about
collections, we will use the language of first order logic 
$L=\{\in\}$ \\ where $\in$ is a binary predicate symbol which has a clear
intuitive meaning: we write ``$x\in y$'' to say that collection $x$ belongs to
the collection $y$. However, as it is well known, it is impossible
to have all the collections in the domain of the language, but only
some (the goal is to have as much as possible), which will be called
{\bf sets}. So the model of the language $L$ will be the (intended)
universe of sets $V$ equipped with the membership relation $\in$. A priori,
it is an arbitrary model of $L$. Each formula $\varphi (x)$ of
the language $L$ is assigned a collection of objects from $V$ satisfying the
formula. The collection will be denoted $\{ x\mid\varphi (x)\}$:

$$ a\in \{ x\mid\varphi (x)\}\ \lra\ V\models \varphi (a)$$

 \noindent  Such collections will be called  {\bf classes}. A
 ``serious'' universe should represent every such  class  by
 its object, set $s$ with the property

$$a\in \{ x\mid\varphi (x)\}\ \lra\ V\models a\in s$$

 \noindent or expressed in the language $L$ itself:

$$V\models \forall x\ (x\in s\ \lra\ \varphi (x))$$

 \noindent  But {\em no matter how we imagine the universe, there
 will be classes of its objects which cannot be represented by its objects}.
An example is {\bf Russell's class} $R=\{ x\mid x\not\in x\}$.
Namely, it is a {\em logical truth} of the language $L$ that there
is no set $R$ with  property

$$\forall x\ \ \ x\in R\ \lra\ x\not\in x$$

 \noindent  Indeed, let us suppose that $R$ is a set. Then, investigating whether it
 is an element of itself, we  get a contradiction:

  $$ R\in R \ \lra\ R\not\in R$$

\noindent  Not only is the result  intuitively unexpected but it is a
{\em logical truth} which we certainly do not understand as well as,
for example, the logical truth $\forall x \ x=x$. The goal of this article is to
 understand better the logic which does not permit  some classes to be
sets.

  The description of a situation  follows. To every set $s$, as the object of the
language $L$, we can associate class $\{ x\mid x\in s\}$.  But, the
opposite is not true. There are classes that can be described in the
metalanguage, like Russell's class, which are not objects of the
language $L$, therefore are not sets. We will call them {\bf proper classes}
or, more in accordance with basic intuition, {\bf paradoxical
classes}.

  To simplify, we will unite the reasoning about sets and
classes into language $LL=\{\in \}$ \label{LL}. This language has the same vocabulary as $L$ but the intended interpretation is different; its objects are  classes of sets. The classes that can be represented with sets  will be identified with sets. Other classes, which cannot be the objects of the language $L$,  we will call {\em paradoxical or proper classes}. Formally, we define in $LL$

\

$A$ is a  {\bf proper} or {\bf paradoxical class} $\leftrightarrow$  $\forall X(A\not\in X)$

\

$A$ is a \textbf{set} $\leftrightarrow$ $A$ is not a proper class.

\

In order to make the translation from the original language $L$ to
the extended language $LL$ easier, we will use capital letters for
variables over classes and small letters for variables over sets.
So, for example, the formula of language $L$ ``$\forall x\varphi
(x)$'' can be considered as an
abbreviation of the formula ``$\forall X (X\ is\ a\ set\ \ra\
\varphi(X))$'' of the language $LL$.

  We consider classes as collections of sets. Every condition on sets
determines the class of all sets satisfying the condition. If
classes have the same members, we consider them equal. These ideas will be formulated in the
 following axioms in the language $LL$ \label{ak}:

\noindent {\bf axiom of extensionality:}

$$A=B\quad\lra\quad \forall x(x\in A\ \lra\ x\in B)$$

\noindent {\bf axiom schema of impredicative comprehension:}

$$\exists A\forall x(x\in A\ \lra\ \varphi (x))$$

\noindent where  $\varphi (x)$ is a formula of the language $LL$
which does not contain $A$ as a free variable.

The consequence  of the axioms is  that  each formula $\varphi
(x)$  can be assigned a unique class   $\{ x\mid \varphi (x)\}$
of all  sets $x$ for which $\varphi (x)$
holds.

  Let's note that to describe classes we will use impredicative comprehension, where $\varphi(x)$ is any formula of $LL$, instead of predicative
  comprehension where $\varphi(x)$ contains only quantification over
  sets. The impredicative comprehension gives us all  classes we  need
  regardless of the way in which we determine them ( for
details see \cite[page \ 119]{fra}). To facilitate  means of expression, we will take a maximalist approach to classes. The basic problem of any set theory
  is not what classes but rather what sets  there are, and our position about classes is irrelevant for this problem. Almost all
  the reasoning (more precisely, all the reasoning which needs only
  predicative comprehension) can be translated either in the language
  $L$ or as the reasoning about the language $L$, so we do not need to
  mention classes in any way. Such approach would  minimize assumptions
  but it would complicate  means of expression. Because of this we
  chose the language $LL$ with the described axioms. The translation
  to $L$ will be done only in some essential situations where it facilitates
   better understanding of the obtained results. Technique of the
  translation is well known --- a discourse about the class
  $C=\{ x\mid\varphi(x)\}$ can be understood as an abbreviation of
  a discourse about the formula $\varphi (x)$. We will need the
  following translations:

\


$s\in\{ x\mid\varphi (x)\} \quad\mapsto\quad\varphi (s)$

\

 $\{ x\mid\varphi (x)\}\textrm{ is a set}\quad\mapsto\quad\exists s\forall x(x\in s\lra\varphi (x))$

 \

  $\{ x\mid\varphi (x)\}\textrm{ is a paradoxical class}\quad\mapsto\quad\neg\exists s\forall x(x\in s\lra\varphi (x))$

  \

With this terminology we only modeled  a necessity for
differentiation between paradoxical classes and sets. Now we will formulate a simple logical criterion to distinguish between sets and paradoxical classes.

  The condition that  there is no set of all objects that satisfy formula $\varphi(x)$ ($\neg\exists s\forall x(x\in s\lra \varphi
(x))$)  can be expressed in a logically equivalent way which shows an elementary logical
mechanism that does
not permit such a set. This logical biconditional we will call \textbf{the productivity principle}. According to that principle {\em paradoxical
classes are exactly those classes for which there is a productive
choice i.e a way to choose for every subset $s$ of the class an object
$c$ of the class which is out of $s$}:

\begin{center}
\includegraphics[height=3cm]{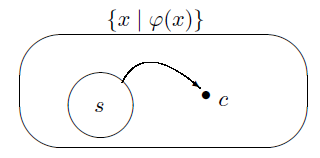}
\end{center}

\

\begin{propozicija} For $\varphi (x)$ a formula of $L$ the following is a logical truth of
$L$:

$\neg\exists s\forall x(x\in s\lra \varphi (x))\ \lra\
  (\forall s(\forall x(x\in s\ra\varphi (x))\ra\exists x(x\not\in s\ and\ \varphi (x))))$
\end{propozicija}

\dzcite ($\ra$): Let $s$ be such that $\forall x(x\in s\ra\varphi
(x))$. According to the assumption on the left side of the biconditional there is an $x$ such that  $x\in s\lra\varphi (x)$
isn't true. From these two conditions it follows that $\varphi (x)\ra x\in s$ is
not true . So $\varphi (x)$ is true  and
$x\not\in s$.

\noindent ($\leftarrow$): Assume the right side of the
biconditional. If there is an $s$ such that $\forall x(x\in s\lra
\varphi (x))$ (*)  then
 $\forall x(x\in s\ra \varphi (x))$. Using the right side of the biconditional, we can conclude that there is an $x$ such that $x\not\in s$ and $\varphi (x)$.
 However,  (*) yields a contradiction  --- $x\not\in s$ and $x\in s$.

\qed

In the language $LL$ this principle has more simple formulation but it
is not a logical truth anymore. Namely, the concept of classes as
collections of basic objects puts a minimal condition
(extensionality of classes) on the principle.

\begin{propozicija}\label{lpo} With the assumption of the axiom of extensionality:

$C$ is a paradoxical class $\lra$ for every set  $s\subseteq C$
there is $x\in C\setminus s$.

\end{propozicija}

\dzcite Suppose that $C$ is a paradoxical class. Then for every set
$s\subseteq C$  \ $s$ is not equal to $C$ (because $C$ is not a set). Hence,  using the axiom of extensionality, we can conclude that  $s$ is a proper subset of $C$ i.e. there is
 $x\in C\setminus s$. Conversely, assume the left side of the biconditional. Suppose that $C$ is
 a set. Then we can take $C$ for a subset $s$ of $C$. Therefore
there is   $x\in C\setminus C=\emptyset$, a contradiction. Hence, $C$ is a paradoxical class.

\qed

A condition for $C$ that for every set  $s\subseteq C$ there is $x\in
C\setminus s$  has been formulated for the first time by Bertrand Russell
\cite{rus} in his generalized contradiction. While Russell's
formulation demands the existence of a definable operation which
 for a given subset of a class gives a new element of the class (a
productive operation on the class), for the productivity principle
the  only thing that matters is the existence of a new object (a productive choice
on the class). Moreover, Russell's view on the meaning of the principle is
different (see page \pageref{russ}).

{\em The productivity principle is a logical truth scheme in the language $L$}. This means that we  assume nothing  about sets and that all of our assumptions about metalanguage $LL$ of classes, introduced to facilitate the discussion, are irrelevant for the principle. The principle tells nothing specific about sets but it  puts  logical bounds on every  theory of sets. {\em The productivity principle is also a simple logical truth}, because we can understand it almost as well as the logical truth $\forall x \ x=x$. Moreover it gives a simple logical mechanism of productive choice which prevents  classes with such mechanism to be sets. A particular theory of sets postulates what sets there are and what operations over sets there are. The productivity
principle logically translates this information into information about what
paradoxical classes there are -  they are collections on which the
postulated (by means of the theory) fund of sets and operations enables a
productive choice.  We will show in the next section how to
produce, using this principle, the known paradoxical classes and find new paradoxical classes. Basically, the principle says that we cannot have  all imaginable sets  and all imaginable operations at the same time. Then we could imagine an operation which when applied to every set gives a new
element and we could also imagine a set that is closed under the operation. This is a
contradiction in itself in the same way as the classical puzzle of the
omnipotence of God. The basic religious intuition is that God is
omnipotent, but then he could make a stone which he could not  move.
The omnipotence requires that he is able to move any stone and that
 he is able to make an unmovable stone. And this is a contradiction. The same conflict of basic intuition
 and logical bounds appears in the naive set theory.

In the sequel we will show how the productivity principle produces paradoxical classes, how it helps  analyze the limitation of size principles and how we can use it as a guide in the construction of a theory about sets.

 \section{Finding paradoxical classes} \label{par}

In this section we  will show how the known paradoxical classes
follow from the productivity principle and how the principle gives us a uniform way
to generate new paradoxical classes
--- {\em we must look for classes having a productive choice on
themselves}. Each paradox (a discovery of paradoxical class) will be
presented in two ways
--- in a uniform way, by establishing a productive choice on the class
(which makes it paradoxical by means of the principle), and in a direct way,
by means of repeating the proof of the productivity principle ( to assume that the
class is a set and to get a contradiction by applying a productive
choice to the class). The first way displays a system to find
paradoxical classes and explicates the logic behind them. The second
way is more direct and explicates a paradoxicality of
such classes with regard   to primary intuition.

 On the logical basis there is only one operation --- identity, $s\subseteq V\ \mapsto\ s\in V$,
so logical paradoxes are associated with classes where identity
is productive:

   $$s\subseteq C\ \mapsto\ s\in C\setminus s$$

   Russell's class is such a class.

\noindent  {\bf Russell's class} $R=\{ x\mid x\not\in x\}$.

  {\em Uniform way}. Let  $s\subseteq R$. If  $s\in s$ then
 $s\in R$, so $s\not\in s$, a contradiction.
Therefore $s\not\in s$. However, then
 $s\in R$. So $s\in R\setminus s$ i.e. identity is productive on $R$.

  {\em Direct way}. Suppose that $R$ is a set. If $R\in R$ then $R\not\in R$.
So, $R\not\in R$. But from the specification of $R$ it follows
that $R\in R$, a contradiction.

The class of not - $n$ - cyclic sets, the class $NI$ of sets which are not isomorphic to one of its elements, and the class of grounded sets are all classes on which identity is productive. Proofs are given in  Appendix. We will establish  here some  relationships between such classes.

Let's note that the universe, in absence of other postulates, is not  a class on which
identity is productive, because we can imagine sets  which contain themselves as members, for example $\Omega =\{\Omega \}\in\Omega$.

  If we write down the productivity condition of identity on a class  $C$ in a more set - theoretical
terminology we have

$$P(C)\subseteq C\cap R$$

\noindent where $R$ is Russell's class. From that it is easy to get
the following results:

\begin{propozicija}
\begin{enumerate}
\item The intersection of classes on which identity is productive is a class on which identity is productive.
\item If we assume the axiom of dependent choices, the class of ungrounded sets $WF$ (for definition see Appendix) is the
      least class on which identity is productive.
\item Every transitive class on which identity is productive  is a subclass of Russell's class.
\end{enumerate}
\end{propozicija}

\dzcite
\begin{enumerate}
\item Let   $P(C_{i})\subseteq C_{i}\cap R$, $i=1,2$. Then
$P(C_{1}\cap C_{2})\subseteq P(C_{i})\subseteq C_{i}\cap R$, so
$P(C_{1}\cap C_{2})\subseteq (C_{1}\cap C_{2})\cap R$.
\item From the groundedness of elements of $WF$, and using the axiom of dependent choices
we can infer the induction principle: if the formula $\varphi (x)$
has the inductive property $\forall x\in y\varphi (x)\
\ra\ \varphi (y)$ then $\forall x\in WF\varphi (x)$. Using this
principle it is easy to prove that every class $C$ on which identity
is productive contains $WF$. Indeed, let  $\forall x\in y\ x\in C$.
Then $y\subseteq C$, so $y\in P(C)\subseteq C$, therefore
 $y\in C$. According to the induction principle we can conclude that every $x$ from $WF$ is in $C$.
\item If $C$ is transitive then $C\subseteq P(C)\subseteq C\cap R\subseteq R$.
\end{enumerate}
\qed

 However, in spite of these results we generally  recognize such classes through their intensional characteristics. All classical examples we have mentioned are like that.

    Other known paradoxes have the origin in the productivity of other operations.
Their list is as follows:

\

 \noindent   {\bf Class $ORD$ of all ordinals}.

  {\em Uniform way}. The productive operation on $Ord$ is to get the first ordinal greater
of all the ordinals from a given set. Namely, from the theory of ordinals
it follows that for every set $s$ of ordinals    such
ordinal exists, let's name it  $s^{+}$. So, $s^{+}\in ORD\setminus s$.

  {\em Direct way}. Suppose  $ORD$ is a set. Then there is the first ordinal out of $ORD$,
 However, it is impossible since, by  definition of $ORD$, $ORD$ contains all the ordinals.

\

 \noindent {\bf Class $CARD$ of all cardinals}.

  {\em Uniform way}.  From the theory of cardinals it follows that for every set of cardinals
 $s$ there is the first cardinal $s^{+}$ greater of all the cardinals from
$s$, so $s^{+}\in CARD\setminus s$.

  {\em Direct way}.  Suppose that  $CARD$ is a set Then there is the first cardinal
$CARD^{+}$ out of $CARD$. However, this is impossible since $CARD$
contains all cardinals.

\

\noindent     {\bf The universe} is paradoxical  if we assume
the subset axiom. Usually, paradoxicality is proved  by means of reduction to
paradoxicality of Russell's class
 $R$. But we can prove it using a suitable productive operation.

  {\em Uniform way}. We will show that diagonalization $s\mapsto \Delta s=\{ x\in s\mid x\not\in x\}$
is a productive operation on the universe (it is enabled by the
subset axiom).
 Indeed, a condition on $\Delta s$ of being an element of itself is

  $$\Delta s\in\Delta s\ \ \lra\ \ \Delta s\in s\ and\ \Delta s\not\in \Delta s$$

\noindent   If $\Delta s\in s$ then we have a contradiction:

   $$\Delta s\in\Delta s\ \lra\ \Delta s\not\in \Delta s$$

  \noindent   Therefore $\Delta s\in V\setminus s$.

  {\em Direct way}. Suppose that the universe  $V$ is a set. Then $\Delta V$
is a set out of  $V$, which is impossible because $V$
contains all sets. Let's note that  $\Delta V=R$ so it is common to
obtain a contradiction by showing that $R$ is not a set.

   We can find another productive choice on the classical universe of sets.
 The partitive set operation enables the choice. By Cantor's theorem there are more
sets in
 $P(s)$ than there are in  $s$, so there is an $x$ in $P(s)$ out of $s$.
Therefore we have a productive choice of $x\in V\setminus s$. If we
transform this argument into a direct proof we obtain another
standard proof for paradoxicality of $V$. Namely, if $V$ is a set
then it  contain $P(V)$ as a subset, but it is impossible
because  $P(V)$ contains by means of Cantor's theorem more elements than $V$.

In a narrow connection with diagonalization and Cantor's theorem is the \v{S}iki\tj 's class \cite{sik}. Its paradoxicality is proved in Appendix.

\

   In the same way for a given fund of operations the productivity principle
gives instructions for generating {\em new} paradoxical classes ---
{\em we need to look for classes on which operations enable
productive choice}.

\

\noindent {\bf class $NWF$ of all ungrounded sets}.

\noindent  If we assume that there is an ungrounded set, that the
universe is a paradoxical class, and the axioms of union and pair,
then we can show that this class is paradoxical.

  {\em Uniform way}. Let $s\subseteq NWF$. If $s=\emptyset$
then we get an element  of $NWF$ (such an element exists by means of  the assumption
of existence of ungrounded sets)
  for a new element of the class. If $s\neq\emptyset$ we get
 $x_{1}\in s$ and $x_{2}\not\in\cup s$ (such an element exists because the universe is not a set).
We claim that $\{ x_{1},x_{2}\}$ is a productive choice on $NWF$.
Really, if $\{ x_{1},x_{2}\} \in s$ then $x_{2}\in \cup s$, contrary
to the choice of $x_{2}$. So $\{ x_{1},x_{2}\}\not\in s$. But
because of ungroundedness of $x_{1}$ there is an ungroundedness of
$\{ x_{1},x_{2}\}\ni x_{1}\ni\ldots$. Therefore $\{ x_{1},x_{2}\}\in
NWF\setminus s$.

  {\em Direct way}. Let $NWF$ be a set. Under the assumption
$NWF\neq\emptyset$ there is $x_{1}\in NWF$. Because $\cup NWF$
isn't equal to the universe there is $x_{2}\not\in \cup NWF$. If $\{
x_{1},x_{2}\}$ belongs to $NWF$ then $x_{2}$ belongs to $\cup NWF$,
contrary to the choice of $x_{2}$. Therefore $\{
x_{1},x_{2}\}\not\in NWF$. But  $x_{1}$ is ungrounded so there is an
ungroundedness of $\{ x_{1},x_{2}\}\ni x_{1}\ni\ldots$,that is to
say $\{ x_{1},x_{2}\}\in NWF$, a contradiction.

        The proof is valid under the weaker assumptions, too. It is enough to suppose, beside the
axiom of pair, that there exists an ungrounded set and that for
every set $s$ of ungrounded sets
 $\cup s\neq V$.

  For the given function $F:V\longrightarrow V$  how can we find a class on which  it
 is productive? If we assume that $F$ is injective, such class is the following

 \

 \noindent {\bf class of all values of injective function  $F$ which do not belong to their
argument} $I =\{ F(x) \mid F(x)\not\in x\}$.

  {\em Uniform way}. We will show that   $F$ itself is productive on $I$.
Let  $s\subseteq I$. If $F(s)\in s$ then $F(s)\in I$. Using the
specification of  $I$  this  means that there is an $x$ such that
$F(s)=F(x)\not\in x$. But the injectivity ensures that  $x$ is equal to $s$,
so we get that $F(s)\not\in s$, and this contradicts  the
assumption. Therefore, $F(s)\not\in s$. However, it implies $F(s)\in I$. Hence,
$F(s)\in I\setminus s$.

  {\em Direct way}. Let  $I$ be a set. If $F(I)$ belongs to $I$, then using the specification of
 $I$ it means that there is  an $x$ such that $F(I)=F(x)\not\in x$.
But using injectivity of  $F$ \ $x$ is equal to $I$, so we get $F(I)\not\in
I$.Therefore, $F(I)$ doesn't belong to $I$. But then, using the
specification of $I$, $F(I)$  belongs to  $I$; a
contradiction.

   Here are some examples of such functions and associated paradoxical classes:

\begin{enumerate}

\item $x\mapsto \{ x\}$:  $I=\{ \{ x\}\mid \{ x\}\not\in x\}$
      (the existence of a singleton is assumed)
\item $x\mapsto \{ x,a\}$:  $I=\{ \{ x,a\}\mid \{ x,a\}\not\in x\}$
      (The existence of a pair is assumed)
\item $x\mapsto x\cup\{ x\}$:  $I=\{ x\cup\{ x\}\mid x\cup\{ x\}\not\in x\}$
      (The existence of  singleton and union is assumed, as well as the foundation axiom)
\item $x\mapsto (x,a)$:  $I=\{ (x,a)\mid (x,a)\not\in x\}$
     (the existence of a pair is assumed)
\item $x\mapsto (x,a,b\ldots)$:  $I=\{ (x,a,b,\ldots )\mid (x,a,b,\ldots )\not\in x\}$
      (the existence of a pair is assumed)
\item $x\mapsto P(x)$:  $I=\{ P(x)\mid P(x)\not\in x\}$
      (the existence of a partitive set is assumed)
\item $x\mapsto F[x]$, for injective  $F$:  $I=\{ F[x]\mid F[x]\not\in x\}$
      (the replacement axiom is assumed)
\end{enumerate}..........

 \

 Paradoxes connected with injective functions can be generalized
by switching from  functions to systems of rules \cite{ach} which are
closer to the idea of a productive choice. The system of rules is a
concept which connects different operations in a unity,
heterogeneous ways of getting new elements into a unique system. Abstractly formulated,
the {\bf system of rules} is every class of ordered pairs.
The intended interpretation is suggested by means of  suitable
terminology. If an ordered pair $(s,x)$ belongs to the system of
rules $\Phi$ we write $\Phi: s\vdash x$ and we say that in  system
$\Phi$ object $x$ is derived from  set of objects $s$. When it is
clear from context what system $\Phi$ is  considered, its label will be dropped.
 Examples of such systems
are formal proof systems, or the system for generating natural
numbers. We usually  use them for generating some objects (let's say
theorems or numbers) using an iteration of rules beginning from some
initial objects (let's say axioms or  number 0). We say that  {\bf
the system is productive on
 $C$} when for all $s\subseteq C$ there is an $x$ such that
$s\vdash x\in C\setminus s$. Although we can generally infer  more
than one object from a given set
 (there is no uniqueness) the concept of injectivity can be easily generalized
to systems --- we say that a system is {\bf deterministic} if every
object can be derived from   one set at the most:

 $$s_{1}\vdash x\ and\ s_{2}\vdash x\ \ \ra\ \  s_{1}=s_{2}$$

\noindent We say that a system is  {\bf global} $\lra\ \forall
s\exists x\ s\vdash x$.

 Let's show that for a global deterministic system  the following class is paradoxical:

 \

\noindent  {\bf class of all objects which do not belong to sets
they are derived from in  a global deterministic system $\Phi$}:

 $I=\{ x\mid \exists a(a\vdash x\ and\ x\not\in a)\}$.

  {\em Uniform way}. Let $s\subseteq I$. Because of globality  there is an $x$ such that
 $s\vdash x$.
If $x$ belongs to $s$ then it  belongs to $I$, so there is an
$a$ such that $a\vdash x\ and\ x\not\in a$. However, using determinism of
the system, $a=s$, so $x\not\in s$, and this is a contradiction.
Therefore $x\not\in s$. But then using the specification of $I$, $x\in
I$. Therefore $x\in I\setminus s$.

  {\em Direct way}.  Let $I$ be a set. Then there is an $x$ such that
$I\vdash x$. If $x$ belongs to $I$ then, using  the specification of $I$
and determinism of the system, $x$ does not belong to $I$, a contradiction. Therefore
$x$ doesn't belong to $I$. But then, using the  specification of $I$, it
belongs to $I$, and this is a contradiction, too.

  Thus, for example, by combining the previous injective operations we can get, with some caution,
new deterministic systems and associated paradoxical classes. For
example, combining operations  $x\mapsto \{ x\}$ and $x\mapsto P(x)$
we get a system which is deterministic because the situation $\{ x\}
=P(y)$ can be obtained only in one way, for $x=y=\emptyset$.

  The concept of productivity can be generalized into the concept of a monotonic operator, too.
For an operator $\phi$ we say that it is  monotonic if for any set
$a$ and  $b$ $\ a\subseteq b\ \ra\ \phi (a)\subseteq \phi (b)$. With
each system  $\Phi$ we can associate a monotonic operator which
maps every set $x$  to the set of all objects inferred  from  some
subset of $x$.

  $$s\mapsto \phi(s)=\{ x\mid \exists a (a\subseteq s\ and\ a\vdash x)\} $$

    Monotonic operator $\phi$ is productive on $C$ if for every
  $s\subseteq C\ \phi (s)\cap C\setminus s\neq\emptyset$. Let's note that
a system needs not to be productive for the associated operator to
be productive. For example the standard system for generating ordinals
is such a system: \label{ORD}

    $$\{\alpha\}\mapsto\alpha\cup\{\alpha\}$$
    $$s\mapsto \cup s$$

  Paradoxicality of $ORD$  doesn't follow from the system because if $\alpha$ is the successor
ordinal then $\cup\alpha\in\alpha$. But in $\alpha$ exists its
predecessor $\beta$, so from $\{\beta\}\subseteq\alpha$ we can infer
$\alpha\not\in\alpha$. By generalization of the previous notice we
can prove that
 the associated monotonic operator is productive and from this  we can
infer  paradoxicality of $ORD$ from it.

    But is there any systematic way to find a class on which a given system of rules is productive,
a class specified by some extensional means and not intensional?
Some results in that direction
 follow.

    Classes on which a given operation $F$ is productive need to be looked for amongst
classes closed under the operation. The first such class is $I(F)$
(so called the least fixed point of $F$), the class of objects given
by iterative application of $F$ beginning with the empty set.
 It can be described in various ways and each of them needs some postulates. We will describe it as
the intersection of all classes closed under $F$:

$$ I(F)=\cap \{ C\mid \forall s\subseteq C \quad F(s)\in C\}$$

 \noindent The description, although intuitively clear, can't be considered literally because elements of
classes can't be classes, but as an abbreviation for \label{pok}

$$ I(F)= \{ x\mid \forall C ((\forall s\subseteq C \ F(s)\in C)\ \ra\  x\in C)\}$$

\noindent In the same way we will consider other similar
descriptions which will be used soon. They are based on the
impredicative principle of comprehension of classes which is
discussed at  the beginning of the article.

   Equally, for  system $\Phi$ we define

$$I(\Phi)=\cap \{ C\mid \forall s\forall x(s\subseteq C\ and\ \Phi :s\vdash x \ \ra\  x\in C)\}$$

 The next proposition shows that the problem of productivity of a given operation $F$
or a given system of rules $\Phi$ on some class is reduced to the
problem of its productivity on its least fixed point: \label{fix}

 \begin{propozicija}

   (i)  For function $F$ there is a class on which it is productive
                        $\lra$ it is productive on $I(F)$.

   (ii) For system of rules   $\Phi $ there is a class on which it is productive
                        $\lra$ it is productive on  $I(\Phi )$

\end{propozicija}

\dzcite: (i): The direction $\leftarrow$ is trivial. Let's prove the
direction $\ra$. Let $F$ be productive on $C$. Then $C$ is closed
under $F$, so $I(F)\subseteq C$. Let $s\subseteq I(F)$. Then
$s\subseteq C$, so, because of productivity on $C$, $F(s)\not\in s$.
However, because
 $I(F)$ is closed under $F$,   $F(s)\in I(F)$, so $F(s)\in I(F)\setminus s$. Therefore, $F$
is productive on $I(F)$. For (ii) the proof is similar.

                             \qed

 \begin{korolar}\label{inj}

 \begin{enumerate}
  \item If $F:V\longrightarrow V$ is injective then $I(F)$ is a paradoxical class.
  \item If $\Phi$ is a global deterministic system of rules  then
         $I(\Phi )$ is a paradoxical class.
 \end{enumerate}
 \end{korolar}

  \dzcite  We have shown  that an injective function $F$ is productive  on
$\{ F(x) \mid F(x)\not\in x\}$ and that a global deterministic
system of rules $\Phi$ is productive on $\{ x\mid \exists a(a\vdash
x\ i\ x\not\in a)\}$. Using the  previous proposition it means that they
are productive on their least fixed point. \qed

For example, the least fixed points of all previously mentioned
injective operations and global deterministic systems are
paradoxical classes.

   \section{Limitation of size principles}\label{toobig}

According to the productivity principle paradoxical classes are exactly those classes on which there is a productive choice. The existence of such a choice doesn't mean necessary that a theory about sets (or the underlying conception of set) assures enough sets and operations. It can even mean the opposite, that there are not enough sets. This is exactly the case with the limitation of size principles. We will see that {\em the limitation of size principles assure a productive choice of a certain kind on paradoxical classes because they are of a restrictive nature and  don't allow more sets}.

 The analysis of paradoxical classes in {\em already existing} Cantor's set theory (which mathematicians
from the beginning of the century as well as today used and
needed) extracted their common property --- they are in  a sense
`` too big''. In Cantor it is in the sense of how many elements
there are in a class (the cardinal sense), in Mirimanoff in the sense of how many stages
there are in  a class formation (the cumulative sense). The idea of ``too big'' is also presented in Zermelo and Fraenkel, although not in such a definite sense.
However, all the variations are in  a harmony with  the
phrase ``too big'' as it is used in ordinary language, so we will use the phrase in a general context. Of
course, if {\em all known} paradoxical classes are ``too big'', it
does not mean that  {\em all} paradoxical classes are as such. Thus the
observation is formulated as  the {\bf  limitation of size
hypothesis} \cite[page \ 176]{hal}: {\em all paradoxical classes
are (in  some sense) too big.} Because the idea of ``too big''
suggests an explanation why some classes are not sets, it is taken
as a  criterion to distinguish between sets and paradoxical classes:

\begin{center} {\em $C$ is a paradoxical class $\lra$ $C$ is too big}  \end{center}

\noindent We will use the word ''small''  instead of "not too big", so in the contraposition form we have

\begin{center} {\em $C$ is a set \ $\lra$ \ $C$ is small}  \end{center}

For every meaning of "too big" we have a corresponding principle, called the {\bf limitation of size principle.}. As every such principle  sounds as an explanation what a set is, it can be taken as a basis of a definite conception of set, which we will call here the {\bf  limitation of size conception}. \footnote{in the  literature on set theory it is common to  use this term in the specific meaning where "too big" is in Cantor's sense} However, without other explanations of what is small (= not too big) the limitation of size principles do not carry any information about sets; we have only new names for the old concepts ( ''small''  for ''set'' and ''too big'' for ''paradoxical class''). Other explanations or postulates must declare what is small. In Cantor, ''small'' means enumerated by some ordinal, and in Mirimanoff ''small'' means to be a subset of some stage $V_{\alpha}$ of the cumulative hierarchy. Of course, what is small depends on what ordinals i.e. what stages  there are. Taken together these other postulates  give us the meaning  of the limitation of size principle in one direction:

\begin{center} {\em $C$ is small  $\ra$ $C$ is a set }  ( postulates on what is small) \end{center}

\noindent It means that the limitation of size principle carries specific information only in the opposite direction, so we will consider the limitation of size principles to be the following statements:

\begin{center} {\em $C$ is too big $\ra$ $C$ is a paradoxical class}  \end{center}

\noindent or in the contraposition form:

\begin{center} {\em $C$ is a set $\ra$ $C$ is small}  \end{center}

We will now analyze the meaning of the limitation of size principles formulated in this way. Other postulates say to us what sets  there are. Acording to the productivity principle these sets and operations on them give us productive choices on some classes and make them paradoxical. Bertrand  Russell was the first to notice this in his correspondence with Jourdain (see  \cite[page \ 180]{hal}). Russell starts from an observation that all known paradoxical classes have a common property which he expressed in
 {\bf Russell's generalized contradiction}\label{russ}:

\

  {\em Let there be  a function symbol $F$  and a formula $\varphi (x)$ such
  that}

  $$\forall s[\forall x(x\in s\ra\varphi (x))\ \ra\ F(s)\not\in s\ and\ \varphi (F(s))]$$

 {\em Then $\{ x\mid \varphi (x)\}$ is a paradoxical class}.

\

  \noindent We can say it in a more readable way:

\

  {\em Let there be a function $F$ and a class  $C$ such that for every
$s\subseteq C\ \ F(s)\in C\setminus s$. Then $C$ is a paradoxical
class}.

\

\noindent  This condition for paradoxicality is the same as in the productivity
 principle, except that in Russell there must be some definable operation that gives
 productive choices. Also, for Russell the condition was not a
 {\em criterion} for
 paradoxicality but a step in arguing limitation of size hypothesis.
 He observed that not only $ORD$ has such a property, but that all
 other classes with such a property contain $ORD$,  or better to
 say, its isomorphic copy. Namely, in a
paradoxical class  the sequence of different elements  can be
reproduced (it is guaranted by  the postulates on what is small):

$$\emptyset,\ F(\emptyset ),\ F(\{ \emptyset ,F(\emptyset )\} ),\ldots$$

 \noindent  The next new element in the sequence is the result of applying $F$ on the set of previous
elements. So we get an ordered structure isomorphic to $ORD$. Even
more, if we take Russell's class and identity we get just
``official'' von Neumann's ordinals. In that way Russell argued
limitation of size hypothesis showing that  all the known
paradoxical classes contain an isomorphic copy of $ORD$, so they are
too big likewise $ORD$ is too big, and not small like sets postulated by  the postulates on what is small. The limitation of size hypothesis saying that all paradoxical classes are too big means in the contraposition form that all small classes are sets. Thus, Russel's argument is the argument that makes plausible the consistency of the postulates on what is small enough to be a set.

Limitation of size principles go further. They say that {\em all classes that are not small are paradoxical classes i.e. there are no more sets except those that are small}. The postulates on what is small say what sets there are and limitation of size principles say that there are no more sets. Thus {\em limitation of size postulates are the closure postulates of the postulates on what is small}. Therefore, they are of a restrictive nature. On  one side, according to the postulates on what is small,
there are collections small enough to be sets. On the other side these postulates assure productive choices on some collections and make them paradoxical and too big. However, there are a lot of collections between these opposites which are not small in the sense of underlying conception on what  small is and which ones are not paradoxical. Let's illustrate this on the cumulative conception. All collections which are built in the cumulative manner beginning with the empty set are sets. The collection which has only itself as  a member, $\Omega =\{ \Omega \}$ doesn't belong to the cumulative hierarchy but it is not too big in the intuitive sense. Furthermore it is not paradoxical. Empty set and $\Omega$ itself are the only subcollections of $\Omega$. If we accept $\Omega$ as a set then $\Omega$ has  two subsets. We can find an object in $\Omega$ outside $\emptyset$ (it is $\Omega$) but not outside $\Omega$. So, there is no  productive choice on $\Omega$.  However, the limitation of size principle will proclaim $\Omega$ too big i.e. paradoxical. According to the productivity principle this means that there is a productive choice on $\Omega$. However, this choice is trivial and enabled by the decision that $\Omega$ is not a set. Namely, because $\Omega$ is not a set the only subset is $\emptyset$ and now we have a productive choice. As another example let's take the whole universe $V$. Although it is not the element of the cumulative hierarchy it is not paradoxical in the absence of other postulates neither. The  consequence of the postulates on what is small is  that the whole cumulative hierarchy $WF$ is paradoxical class, as it is well known (see  Appendix). If we allow some collections which are outside the cumulative hierarchy, like $\Omega$, to be sets,  then $V$ is not necessarily paradoxical, but if we don't allow more sets, and this is the meaning of the limitation of size principle, then $V=Wf$  is a paradoxical class.

On the previous examples we can see that the cumulative limitation of size principle enables a productive choice because it doesn't allow some collections to be sets. Now we will show that this is generally so. We will show that  the limitation of size principle  enables a productive choice of a certain kind on paradoxical classes. Basically, this  choice is  a choice of an object $x$ outside  $V_{\alpha}$,  and  it is a productive choice because the limitation of size principle doesn't allow  collections outside the cumulative hierarchy  to be sets.

Let's suppose some postulates about cumulative hierarchy (the postulates on what is small). The exact formulation of these postulates doesn't matter here. It is enough that they say that every subcollection of some stage $V_{\alpha}$ is a set. Then, from the cumulative principle of the limitation of size ("too big" means now not to be a subset of a some stage $V_{\alpha}$)

$$\forall\alpha\exists x \ x\in C\setminus V_{\alpha} \ \rightarrow \  C\textrm{  is a paradoxical class}$$

\noindent it follows the statement about the existence of a certain kind of a productive choice on a paradoxical class:

$$C\textrm{  is a paradoxical class }\lra \forall s\subseteq C\  \exists  \alpha\exists x \  s\subseteq V_{\alpha}  \textrm{ and }  x\in C\setminus V_{\alpha}\textrm{ (*)}$$

Proof. Let's suppose the limitation of size principle. The direction $\leftarrow$ of (*) follows from logic. Namely, according to the productivity principle we must prove that there is a productive choice on $C$ and the left side of (*) provides such a choice. To prove the direction $\ra$ let's suppose that $C$ is a paradoxical class and that $s\subseteq C$. From the limitation of size principle we can conclude by contraposition that there exists $\alpha$ such that $s\subseteq V_{\alpha}$ (there are no other sets!). However, $C$ is a paradoxical class, so, according to the postulates about cumulative hierarchy, it is not a subcollection  of $V_{\alpha}$. Thus there is $x\in C\setminus V_{\alpha}$. Therefore, we have proved the direction $\ra$ of (*).

The same analysis can be carried out for the Cantor's cardinal limitation of size principle. The principle proclaims paradoxical all collections that cannot be enumerated by an ordinal (= too big collections):

\

\noindent {\em There is no ordinal $\alpha$ and function $F$ such that $C=F[\alpha ]$ \\ $\ra$ $C$ is a paradoxical class}

  \

\noindent Again, we have  postulates about what is small (postulates about ordinals, the replacement postulate, etc.). We will call them  postulates about enumerated collections. These postulates assure  productive choices on some collections which makes them paradoxical and not enumerated by ordinal ($ORD$, the universe $V$, etc.). However, what is with the collection $V_{\omega}$ of all hereditary finite sets, for example? If we don't postulate the existence of infinite sets then $V_{\omega}$ is not enumerated, but there is  nothing which makes it paradoxical. For any finite subcollection (therefore a set) $s$ we can find an element of $V_{\omega}$ outside $s$ (because $s$ is finite and $V_{\omega}$ is infinite), but for the whole $V_{\omega}$ we cannot. If we close the conception of hereditary finite sets by corresponding cardinal limitation of size principle which says that there are no other sets then $V_{\omega}$ is not a set. The finite subcollections of $V_{\omega}$ are the only subsets and now we have a productive choice. Such a choice is again the consequence of the decision not to allow more sets. Again, we can prove that from the cardinal limitation of size principle (from forbidding some collections to be sets) follows the existence of a productive choice of a certain kind on paradoxical classes:

$$ C\textrm{ is a paradoxical class }\lra \ \forall s\subseteq C \ \exists F\exists\alpha \ \ F[\alpha ]=s \textrm{ and } F(\alpha ^{+})\in C\setminus s$$

We can repeat the same analysis for Zermelo's approach to paradoxes. In his approach "small" means to be a subcollection of a set (the subset axiom) and the relative character of the meaning of "small" is explicit here. The corresponding limitation of size principle is

\begin{center} {\em there is not $s$ such that $C\subseteq s$ \ $\ra$ \ $C$ is a paradoxical class} \end{center}

\noindent However, the principle is now the truth of logic, therefore it says nothing about sets. Namely, in a contraposition form  we need to prove that if $C$ is a
set, then there is a set $s$ such that $C\subseteq s$. And  the
assumption gives such a set - it is just $C$. The presence of a productive choice on a paradoxical class follows directly, too.
Because such a class $C$ goes beyond every set, it goes also beyond
any of its subsets, so for $s\subseteq C$ there is   $x\in
C\setminus s$.

Using the productivity principle we can analyze Fraenkel's limitation of size explanations of the concept of set which are rather loose.   Explaining the difference between paradoxical classes and sets Fraenkel \cite{fra1,fra2}
says that paradoxical classes are of {\em unbounded extension} and
also that they {\em involve too much}.  For the first metaphor of
an unlimited extension we  consider it similar to Zermelo's
metaphor which we have already analyzed. Concerning the second metaphor of ``too big
involvement'', as the main argument Fraenkel mentions that sets, as
opposite to paradoxical classes, can be extended (the idea is
present already in Cantor) and he illustrates this with
diagonalization (which requires the subset axiom):

$$s\mapsto \Delta s=\{ x\in s\mid x\not\in x\}\not\in s$$

 \noindent Conversely to sets, neither Russell's class nor the universe are extended under diagonalization.
For example, for $s\subseteq R \ \ \Delta s=s\not\in s$, so $\Delta
s\in R$.
 It is not clear from Fraenkel's
exposition whether he considers the closure under diagonalization to be a
distinguishing property, or the closure under any operation which
extends sets. The first thought  is wrong because for example $ORD$
is a paradoxical class which is not closed under diagonalization
 (for $s\subseteq ORD$ which is not an ordinal $\Delta s=s$ doesn't belong to $ORD$).
The second thought has the confirmation in Fraenkel's metaphor of
the ``wall'' according to which,
 contrary to sets, to constitute Russell's class ``elements have to be taken from
outside every wall, no matter how inclusive''(\cite{fra1}). This gives
a distinguishing criterion which is just the productivity principle
and Fraenkel's explanation of the concept of set is a logical truth. So, it says nothing specific about sets.

 Von Neumann's  limitation of size principle is a variation of the cardinal limitation of size principle. In this approach "too big" means equipotent to the universe $V$ of all sets:

\begin{center}{\em $C$ is equipotent to the universe $V$ \ $\lra$  \ $C$ is a paradoxical class} \end{center}

\noindent The postulates on what is small will tell what collection are  not equipotent to the universe and hence are sets. Because $ORD$ is a paradoxical class,  it is equipotent to the universe $V$, and as such, by the postulates on what is small, to  every paradoxical class. It means that all paradoxical classes have the same kind of a productive choice as $ORD$ has.

 \section{Heuristic of the productivity principle}

The productivity principle can guide us in the construction of the consistent set theory in the following way. When some sets and operations on sets are postulated the  principle tells us what classes are paradoxical. We cannot proclaim any such class to be a set because we would get a contradiction. However, we can proclaim some other classes, which are not paradoxical by means of the principle, to be sets. So we get a richer theory. Now we can repeat the process and get more sets, and so on. If we do that in some systematic manner we can get a pretty rich theory. Of course,  we cannot prove that  such a theory is consistent but the procedure would convince us that it is the case. We will sketch  one such a theory, called here \textbf{cumulative cardinal theory of sets} which realizes the corresponding cumulative cardinal conception of sets. Because of the construction its consistency seems plausible. Moreover, it inherits good properties from cardinal conception and cumulative conception of sets and as such it easily  justifies  ZFC axioms. Because of the cardinality principle it justifies easily the replacement axiom, and because of the cumulative property it justifies easily the power set axiom and the union axiom.

Before the construction we need to  repeat (page \pageref{fix}) that from the productivity of a function $F$ follows
that $I(F)$ is a paradoxical class. From the
assumption that
 $I(F)$ is a set and applying $F$ on $I(F)$ we get a contradiction. For the  productivity of
 $F$ it is necessary that $F$ is defined on every subset of
$I(F)$. But some functions are such that it is not always fulfilled
--- there are given conditions of their applicability. If $I(F)$ belongs to the
condition, although the function  isn't defined on every subset of
$I(F)$, we will again get a contradiction and therefore show that $I(F)$
is a proper class.
 Contrary, if $I(F)$ doesn't satisfy
the condition the argument fails. Moreover, it seems plausible just
the opposite
 --- that  $I(F)$ is a set. Namely, $I(F)$ is a collection generated by iterative application of $F$,
 so if the operation itself doesn't give a contradiction it would be surprising that some other
operation or system of rules does so.

We will  collect  objects into
sets  in an order of
growing cumulative cardinality of collecting. The collecting of the previous kind produces the least fixed point which will be
considered as an initial place for the collecting of the next kind which has a greater cardinality. The basic
idea is the following. Beginning with the empty set we can get all
hereditary finite sets by finite and combinatory means, so their
existence is indisputable. For example, we can generate them with the
following system of rules:

 \

$ \emptyset\vdash \emptyset$\\ \indent (we accept the empty set)

 \

$\{ x,y\}\vdash x\cup \{ y\}$ \\ \indent (by adding to set $x$ one
element $y$ we get again a set)

 \

\noindent  All sets  generated in this way are hereditary finite (they are not
only finite, but their elements are finite too, and elements of
their elements too, etc.) and they make the least fixed point
$C_{1}$ of the operation of  {\em finite} collecting of objects to
sets. However, $C_{1}$ is not finite, so the operation isn't applicable
to it. Using the heuristic of the productivity principle we accept $C_{1}$
as a set. However, accepting it as a set we accept a countable collecting of
objects into sets, too:

\

$X\subseteq V\vdash X\in V \qquad\textrm{when there is an $F$ such that $F[C_{1}]=X$}$

\noindent where $F[S]=\{ F(x)\mid x\in S\}$, and it need not to be
that $S\subseteq dom(F)$. If $\exists F\ F[Y]=X$ we will write $X\preccurlyeq Y$ ($X$ is dominated by $Y$ or $X$ is $Y$-enumerated)

\noindent  This type of argument that accepting a set with some
cardinality means  accepting all collecting with that cardinality we
will call the {\bf cardinality principle}.

  Beginning with objects from $C_{1}$ and iterating
countable collectings ( it includes finite collectings too) we get
all hereditary countable sets. They make the least fixed point
$C_{2}$ of the collecting. $C_{2}$ contains $C_{1}$ as a subcollection.
But because $C_{1}\not\in C_{1}$ and $C_{1}\in C_{2}$, \ $C_{2}$ is a richer
collection:

                   $$C_{1}\subset C_{2}$$

  The story repeats itself now. Let's show using a modification  of {\em Cantor's argument}\label{cd}
that a countable collecting isn't applicable on
 $C_{2}$. If it isn't  so  then there is the function  $F$ such that
$F[C_{1}]=C_{2}$. Let's consider $D=\{ x\in C_{1}\mid x\not\in
F(x)\}$. Elements of $D$ are sets from $C_{1}\subseteq C_{2}$ and
$D=G[C_{1}]$, where $G$ is identity function on $D$, so  $D$ is a
set from $C_{2}$. However, using the assumption on $F$, there is
$d\in C_{1}$ such that $F(d)=D$. This leads to a contradiction ---
$d\in F(d)\ \lra\ d\not\in F(d)$. Therefore, a countable collecting
isn't applicable on the class of all hereditary countable
 sets $C_{2}$.

  Because a countable collecting isn't applicable on class $C_{2}$ of all hereditary countable
sets,  using the heuristic of the productivity principle we can accept
it as a set. However, $C_{2}$ has a greater cardinality, so accepting it
as a set, using the cardinality principle we accept all collectings of
that cardinality and the process goes on.

  In the same way we solve every next finite jump because the previous proof remains valid
when $C_{1}$ is substituted by $C_{\alpha}$, $C_{2}$ by $C_{\alpha ^{+}}$ (where $\alpha ^{+}$ is the next ordinal following $\alpha$), and the expression ``countable'' by ``{\bf
$C_{\alpha}$-countable}'' (enumeration by $C_{\alpha}$).

 Accepting $C_{\alpha}$ for a set, using the cardinality principle we accept  any
$C_{\alpha}$-collecting of {\em already collected} objects, too. By means of the
previous considerations $C_{\alpha}$-collecting is not applicable on its fixed point
$C_{\alpha ^{+}}$ and using the heuristic of the productivity principle
we accept $C_{\alpha ^{+}}$ as
 a set. By repeating this process we get  an ascending hierarchy of sets:

$$\emptyset = C_{0}\subset C_{1}\subset C_{2}\subset\ldots
             \subset C_{\alpha}\subset C_{\alpha ^{+}}\ldots$$

  A transfinite jump is obtained by  putting together all of the already accepted $C_{\beta}$-collectings,
 over all ordinals $\beta$ less than the limit ordinal $\alpha$. Such a collecting will be called
 {\bf $<C_{\alpha}$-collecting}:

$$X\subseteq V\vdash X\in V \qquad \textrm{when there is $\beta <\alpha$, such that $X$ is $C_{\beta}$-enumerated}$$

 \noindent The least fixed point of the collecting gives the next limit member $C_{\alpha}$ of the hierarchy.
  From the definition of limit $C_{\alpha}$  it follows that $C_{\alpha}$ contains all previous stages and
that it is not identical to any of them (because when it contains
some stage $C_{\beta}$
 it contains also the next stage which $C_{\beta}$ doesn't contain).
Let's show that  $C_{\alpha}$ isn't $<C_{\alpha}$-enumerated. If it
is so, then there is the function $F$ such that
$F[C_{\beta}]=C_{\alpha}$ for some $\beta <\alpha$. As we did for a finite
jump, let's consider $D=\{ x\in C_{\beta}\mid x\not\in F(x)\}$.
Elements of $D$ belong to $C_{\beta}\subset C_{\alpha}$ and $D$ is
$C_{\beta}$-enumerated, so it belongs to $C_{\alpha}$.
 Therefore, there is  $d\in C_{\beta}$ such that $D=F(d)$.
However, this entails the contradiction:  $d\in F(d)\ \lra\
d\not\in F(d)$.

  Because $<C_{\alpha}$-collecting  isn't applicable on the limit  $C_{\alpha}$ \  and $C_{\alpha}$
 is its least fixed point, using the heuristic of the productivity principle we accept $C_{\alpha}$ as a set.
So we can continue the sequence through transfinite jumps, too:

$$\emptyset = C_{0}\subset C_{1}\subset C_{2}\subset\ldots
\subset C_{\alpha}\subset C_{\alpha ^{+}}\ldots\subset
C_{\omega}\subset\ldots$$

\noindent Every stage of the hierarchy contains greater
cardinalities than previous stages and it is a measure for yet
greater cardinalities. The hierarchy spreads as much ordinals as we
can imagine. But as the construction gives greater cardinalities, at the same time  it
gives ordinals for next jumps.

Let's note some properties of the stage $C_{\alpha}$. Because $C_{\alpha}$ is the least fixed point of  collectings with smaller cardinality there is a corresponding induction principle for proving properties of $C_{\alpha}$ and a corresponding recursion principle for defining functions on $C_{\alpha}$. Furthermore, the following sentences  are true:

$$S\subseteq C_{\alpha }\ \ra\ S\in C_{\alpha ^{+}}$$

$$S\in C_{\alpha }\ \ra\ S\subseteq C_{\alpha }$$

$$S\in C_{\alpha ^{+}}\quad \lra\quad S\subseteq C_{\alpha ^{+}}\ and\ S\preccurlyeq C_{\alpha }$$

\

\noindent For limit ordinal $\alpha$

\begin{center} $S\in C_{\alpha }\quad \lra\quad S\subseteq C_{\alpha }\ and\ S\preccurlyeq C_{\beta }$ for some $\beta <\alpha$ \end{center}

Can we go further, by putting together {\em all} $C_{\alpha}$-collectings for the final jump:

$$X\subseteq V\vdash X\in V \qquad \textrm{when there is an $\alpha$, such that $X$ is $C_{\alpha}$-enumerated}?$$

\noindent If we could get new sets in this way   then, by the same argument as in previous jumps, the universe $V$ of all sets would be just a new stage for new jumps. This is impossible, by the very idea of the universe of all sets. Therefore the hierarchy of $C_{\alpha }$ gives all sets and it is closed under $C_{\alpha}$-collectings:

$$V=\bigcup_{\alpha\in ORD}C_{\alpha}=\{S|\exists C_{\beta }(S\preccurlyeq C_{\beta })\}  $$

\noindent It means that

$$ S\textrm{ is a set }\lra \ \exists C_{\alpha }(S\preccurlyeq C_{\alpha })\ \lra\ \exists C_{\beta }(S\in C_{\beta })  \lra\ \exists C_{\gamma}(S\subseteq C_{\gamma })$$

This gives  the limitation of size principle for the cumulative cardinal conception of set:

\

 {\em If $C$ is not a subcollection  of some stage $C_{\alpha}$ (=not dominated by some stage $C_{\beta}$) then $C$ is a paradoxical class}

 \

 Of course, this limitation of size principle is a closuring principle for this conception of set as well as other limitation of size principles are. For example, it excludes non- well- founded sets.\footnote{Let's note that we can replace the least fixed points in the cumulative cardinal hierarchy with the greatest fixed points and get non-well-founded sets, too.}

 We have described the {\bf cumulative - cardinal conception of  set} informally. Now, we will do it more formally. We will formulate axioms in the language $LL$ (page \pageref{LL}) and call them the {\bf CC axioms}. For classes we accept the {\bf  axiom of extensionality} and the {\bf  axiom schema of impredicative comprehension} (page \pageref{ak}). Also, we accept the {\bf axiom of choice}:
\

{\em For every set $s=\{ x_{i}|i\in a\}$ of nonempty disjoint  sets there is a set $c$, called a choice set for $s$,  such that for every $i\in a$ \ $x_{i}\cap c $ contains exactly one element}.

\

\noindent The cumulative cardinal conception is neutral about the axiom of choice as well as cardinal and cumulative conceptions, because the axiom of choice is of a different nature.  However, the cumulative cardinal conception permits the axiom of choice. By the axioms of union and subset (which are valid in the cumulative cardinal conception, as we shall see), if the collection $c$ exists, then $c\subseteq  \cup\{ x_{i}|i\in a\}$ and $c$ is a set.

Before describing the cumulative - cardinal hierarchy of sets  we must describe {\bf ordinals} ---
supports of the construction. For that purpose we need the following elementary axioms for sets:

{\bf  axioms for hereditary finite sets}

   \begin{enumerate}
   \item $\emptyset$ is a set.
   \item If $s$ and $a$ are  sets then $s\cup\{ a\}$ is a set.
   \end{enumerate}

\noindent  These axioms are obviously true in the cumulative - cardinal conception. Now we can define

\

  $s^{+}=s\cup\{ s\}$

\

\noindent and we can define the  class $Ord$ (see page \pageref{ORD})

\

$Ord = \cap\{ C\mid\emptyset\in C\ and\ \forall\alpha (\alpha\in C\ \ra\ \alpha ^{+}\in C)\ and\
 \forall s\subseteq C(s =\cup s\ \ra\ s\in C)\}$

\

 \noindent Intuitively, ordinals grow up together with the cumulative hierarchy. Formally, the next axiom,
 the axiom of cumulative hierarchy,  gives all classical ordinals.

\

the {\bf axiom of a cumulative - cardinal hierarchy}:

\

There is $C:Ord\longrightarrow V$ such that:

\begin{enumerate}
\item $C_{0}=\emptyset$
\item $C_{1}=\cap\{ Cl\mid\emptyset\in Cl\textrm{ and }\forall s,a(s,a\in Cl\ \ra\ s\cup\{ a\}\in Cl)\}$
\item $C_{\alpha ^{+}}=\cap\{ Cl\mid\forall S\subseteq Cl(S\preceq C_{\alpha }\ \ra\ S\in Cl)\}$, for $\alpha > 1$
\item For limit ordinal $\alpha$ $C_{\alpha }=\cap\{ Cl\mid\forall S\subseteq Cl(S\preceq  C_{\beta }
\textrm{ for some $\beta <\alpha$} \ \ra\ S\in Cl)\}$
\item $V=\bigcup_{\alpha\in ORD}C_{\alpha}=\{ S|\exists C_{\beta }(S\preccurlyeq C_{\beta })\}  $
\end{enumerate}

\noindent  The meaning of the relation  $\preceq $  is standard: $X\preceq  Y\ \lra\ \exists F\ F[Y]=X$

Developing the set theory from these axioms is out of the scope of this article. Therefore, we will return to the informal  cumulative - cardinal conception of  set, based on basic ideas of the cardinality of collecting and the possibility to
iterate  greater and greater collectings.  The conception inherits good properties from cardinal conception and
from cumulative conception of sets. It is well known that Cantor's limitation of size conception justifies easily the replacement axiom but it cannot justify the power set axiom and the union axiom. Contrary, the iterative conception justifies  the power set axiom and the union axiom easily but it cannot justify the replacement axiom \cite[page \ 199]{hal}\cite[page \ 208]{bar}.\footnote{ For an elaborate discussion of the cumulative and cardinal conceptions of set see \cite{bool1}\cite{bool2}}  Now we can justify  all $ZF$ axioms easily (we accept the {\bf axiom of extensionality} as a basic property of sets and the {\bf axiom of choice} as an axiom of a different nature which is compatible with this conception):

\begin{enumerate}
\item the {\bf  axiom of empty set}: $\emptyset$ is a set.

   $\emptyset \subseteq C_{1} \ra \emptyset \in V$

\item the {\bf  axiom of infinity}:
$\omega  \subseteq C_{1}\ra \omega\in V$.

\item the {\bf  subset axiom}: $C\subseteq s\ \ra\ C$ is a set,

 Because $s$ is a set it is dominated by some $ C_{\alpha}$. But then  its subclass $C$ is dominated
by $ C_{\alpha}$ too, so $C$ is also a set.

\item the {\bf  power set axiom}:\ $s$ is a set $\ra$ $P(s)$ is a set.

Because $s$ is a set $s\subseteq C_{\alpha}$. But then for $a\subseteq s$
 $a\subseteq C_{\alpha}$, so $a\in C_{\alpha ^{+}}$. It means that $P(s)\subseteq C_{\alpha ^{+}}$. Therefore $P(s)$ is a set.

 \footnote{Because the hierarchy $C_{\alpha}$ gives a cardinal scale we can estimate  where  $P(C_{\alpha})$ is in the scale: $C_{\alpha}\prec P(C_{\alpha}) \preccurlyeq C_{\alpha ^{+}}$. In the cumulative cardinal conception it is natural to postulate that the hierarchy $C_{\alpha}$ provides all cardinalities. From that postulate follows the Generalized Continuum Hypothesis}

\item the {\bf  replacement axiom}: $s$ is a set and $F$ is a function $\ra$ $F[s]$ is a set.

  Indeed, $F[s]\preccurlyeq s\preccurlyeq C_{\alpha}$, so $F[s]\preccurlyeq C_{\alpha}$. Therefore $F[s]$ is a set.

\item the {\bf  union axiom} $s$ is a set  $\ra$ $\cup s$ is a set.

Because $s\subseteq C_{\alpha}$, every $x$ from  $s$ is in $C_{\alpha}$. But then
 $x\subseteq C_{\alpha}$ and $\cup a=\cup\{ x\mid x\in a\}\subseteq C_{\alpha}$. Therefore $\cup s$ is a set.

\item the {\bf  axiom of groundedness} Every set is grounded, that is to say there is no infinite
descending sequence of sets $s\ni s_{1}\ni s_{2}\ni\ldots$.

 This is easy to prove by the $\in$-induction on sets (the principle  is valid because every $C_{\alpha}$ is the least fixed point of the collecting of a certain kind. Let every element of $s$ be grounded. If $s$
isn't grounded then there is an infinite sequence $s\ni s_{1}\ni s_{2}\ni\ldots$. Then the sequence
 $s_{1}\ni s_{2}\ni\ldots$ is infinite, too, which is impossible because by the induction hypothesis $s_{1}$
is grounded.

\end{enumerate}

We have got all the axioms of the Morse -- Kelley set theory. Contrary, from the Morse -- Kelley axioms it would be possible to  prove all $CC$ axioms (We can show that $C_{\alpha ^{+}} = H(|C_{\alpha}|)$, where $|s|$ denotes the cardinality of the set $s$, and $H(\kappa)$ denotes the set of all sets of the hereditarily cardinality less or equal to $\kappa$). In this way the CC axioms provide a natural and plausibly consistent axiomatization for the Morse -- Kelley set theory. 
 \section{Appendix}

In the appendix the paradoxicality of some famous classes is shown in a uniform and in a direct way (see page \pageref{par})

\noindent {\bf The class of not - $n$ - cyclic sets} $R^{n}=\{ x\mid \neg x\in ^{n}x\}$
(where  $x\in ^{n}y\ \lra\ \exists x_{1}, x_{2},\ldots ,x_{n-1}\ \
x\in x_{1}\in \ldots \in x_{n-1}\in y$).

  {\em Uniform way}. Let  $s\subseteq R^{n}$. If $s\in s$
then  $s\in R^{n}$, so there is no  $n$-cycle $s\in x_{1}\in \ldots
\in x_{n-1}\in s$. But such a cycle is just $s\in s\in s\in
\ldots\in s$ ($n$ times). Therefore, $s\not\in s$. If  $s$ doesn't
belong to the class $R^{n}$ then there is $n$-cycle $s\in
x_{1}\in \ldots \in x_{n-1}\in s$. But then there is also
 $n$-cycle $x_{n-1}\in s\in x_{1}\in \ldots \in x_{n-1}$, and it is impossible
because $x_{n-1}$ is an element of  $s$, and so of $R^{n}$,
therefore not-$n$-cyclic one. So the final conclusion is that $s\in
R^{n}\setminus s$, that is to say  identity is productive on $R$.

  {\em Direct way}  Suppose that $R^{n}$ is a set. If $R^{n}\in R^{n}$ then it is on one side
not-$n$-cyclic, because it belongs to $R^{n}$, and on the other side
it is $n$-cyclic because it makes $n$-cycle $R^{n}\in
R^{n}\in\ldots\in R^{n}$ ($n$ times). Therefore $R^{n}\not\in
R^{n}$. So, it belongs to some $n$-cycle $R^{n}\in x_{1}\in
\ldots \in x_{n-1}\in R^{n}$ which gives $n$-cycle $x_{n-1}\in
R^{n}\in x_{1}\in \ldots \in x_{n-1}$. Therefore
 $x_{n-1}\in R^{n}$ is $n$-cyclic, and this is a contradiction.

\

 \noindent {\bf Class $NI$ of sets which are not isomorphic to its element}.

 \noindent We say that
$x$ is isomorphic to  $y$ $\lra$  $\exists F:Tr(x)\longra Tr(y)$
where $F$ is a bijection and preserves the belonging, that is  $a\in b\
\lra\ F(a)\in F(b)$, and  where $Tr(x)=\cap\{ C\mid C\textrm{ is a
transitive and } x\in C\} = \{ y\mid\forall C(C\textrm{ is a
transitive and } x\in C \ \ra\ y\in C\}$.

{\em Uniform way}. Let  $s\subseteq NI$. If $s\in s$ then
 $s$ is isomorphic to its element (to itself). However, this is impossible because
$s$ as a subset of $NI$ contains sets which are not isomorphic to
its element. Therefore, $s\not\in s$.  If $s$ is not in $NI$ then it
is isomorphic to its element $x$ by an isomorphism $F$. From $x\in
s$ results that $F(x)\in F(s)=x$, so
 $x$ is isomorphic to its element $F(x)$. But this is impossible because
$x$ as an element of $s$ is also an element of
 $NI$. So $s$ is not isomorphic to its element. Therefore $s\in NI\setminus s$.

  {\em Direct way}. Suppose that $NI$ is a set. If $NI\in NI$ then it is isomorphic to its element (to itself),
so $NI\not\in NI$. But if it is not an element of $NI$ then it is
isomorphic to its element $x\in NI$, which is  isomorphic to its
element, and this is impossible because $x$ belongs to $NI$. Therefore
$NI\in NI$, and this is a contradiction.

\

 \noindent {\bf Class $WF$ of grounded sets }.

 \noindent We say that  set $x$ is grounded $\lra$
there isn't a sequence of sets $x_{n},\ n\in\omega$ such that
$x=x_{0}\ni x_{1}\ni x_{2}\ni\ldots$---   the existence of
natural numbers = finite ordinals is assumed.

  {\em Uniform way}. Let  $s\subseteq WF$. If $s\in s$ then
 $s\in WF$. This is impossible because there is a sequence
$s\ni s\ni s\ni\ldots$. Therefore, $s\not\in s$. If $s\not\in WF$
then there is a sequence
 $s\ni x_{1}\ni x_{2}\ni\ldots$ and so there is also a sequence $x_{1}\ni
x_{2}\ni\ldots$ from which there follows that $x_{1}$ is ungrounded.
However, this is impossible because $x_{1}$ is an element of  $s$, and so of
$WF$. Therefore, $s\in WF\setminus s$.

  {\em Direct way}. Let it $WF\in WF$. Then it is grounded,
but we have a witness of its ungroundedness --- a sequence $WF\ni
WF\ni WF\ni\ldots$. So $WF$ is ungrounded. But then there is an
infinite sequence $WF\ni x_{1}\ni x_{2}\ni\ldots$, and so there is
also an infinite sequence $x_{1}\ni x_{2}\ldots$ which means that
$x_{1}$ is ungrounded. But this can't be because $x_{1}\in WF$.

\noindent {\bf \v{S}iki\tj 's class}  $S=\{x \mid x\not\in F(x)\}$,
for the surjection $F$ on the universe \cite{sik}:

  {\em Uniform way}. For $s\subseteq S$ there is a $d$  such that $s=F(d)$
(by surjectivity of $F$). We will show that  $d$ itself  is a new
element of $S$. If  $d\in s=F(d)$ then $d$ belongs to  $S$ from
which it follows that
 $d\not\in F(d)$. So, $d\not\in F(d)$. But then
$d\in S$. Therefore, $d\in S\setminus s$.

  {\em Direct way}. Suppose that $S$ is a set. Then, by surjectivity of $F$,
there is  a $d$ such that $S=F(d)$. But then a condition for the
belonging of $d$ to set $S$ is $d\in F(d)\ \lra\ d\not\in F(d)$, and
this is a contradiction.

   Examples of such operations are $x\mapsto \cup x$, $x\mapsto \cup\cup x,\ldots$
because for every set  $s$ \ $s=\cup\{ s\} =\cup\cup\{\{ s\}\}
,\ldots$ (the union and the singleton axioms
 are assumed), and also
$x\mapsto \cap x$, $x\mapsto \cap\cap x$ because for every set  $s$
\ $s=\cap\{ s\} =\cap\cap\{\{ s\}\} ,\ldots$ (the intersection and
the singleton  axioms are assumed)\cite{sik}. 

\bibliography{setparadox}


 \end{document}